\documentclass[10pt,conference]{./IEEEtran}
\usepackage{amsmath,amssymb} 
\usepackage{graphicx}
\usepackage{bm}
\newtheorem{theorem}{Theorem}
\newtheorem{definition}{Definition}

\newtheorem{example}{Example}

\newcommand{\qed}{\hfill $\square$ \par}

\newcommand{\rme}{\mathrm{e}}

\newcommand{\llangle}{\langle\kern -.23em \langle}
\newcommand{\rrangle}{\rangle\kern -.23em \rangle}

\renewcommand{\baselinestretch}{0.95}

\begin{document}
\title{Upper Bounds for the Number of Solutions \\ to Spatially Coupled Sudokus}
\author{
  \IEEEauthorblockN{Tadahiro Kitazono}
  \IEEEauthorblockA{Graduate School of Information Sciences, \\
  Hiroshima City University, \\
  Hiroshima 731-3194, Japan\\
  Email:kitazono@cm.inf.hiroshima-cu.ac.jp}
  \and
  \IEEEauthorblockN{Kazushi Mimura}
  \IEEEauthorblockA{Graduate School of Information Sciences, \\
  Hiroshima City University, \\
  Hiroshima 731-3194, Japan\\
  Email:mimura@hiroshima-cu.ac.jp}
}
\maketitle
\begin{abstract}
Based on combinatorics, we evaluate the upper bounds for the number of solutions to spatially coupled Sudokus, 
which are popular logic puzzles. 
\end{abstract}
\IEEEpeerreviewmaketitle


\section{Introduction}
\par
\begin{figure}[t]
  \begin{center}
    \includegraphics[width=40mm]{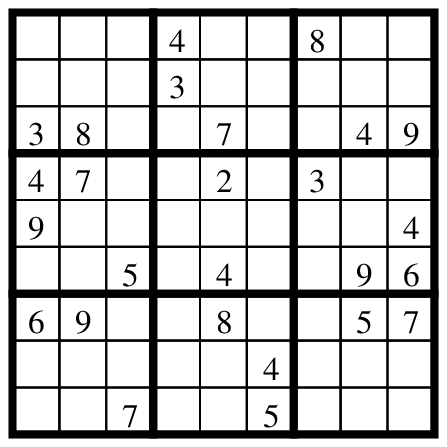}
    \caption{A $3 \times 3$ Sudoku puzzle}
    \label{fig:example}
  \end{center}
\end{figure}

Sudoku is a popular logic puzzle. 
It is presented with a $9 \times 9$ grid, in which some cells have a digit from 1 to 9. 
The task is to complete the grid, by filling in the remaining cells such that each row, each column, 
and each one of nine $3 \times 3$ blocks contains the digits from 1 to 9 exactly once. 
An example of $3 \times 3$ Sudoku is shown in Figure \ref{fig:example}. 
A pattern that fills in the empty cells in a Sudoku 
according to the three constraints mentioned above is called a {\it solution}. 
So far the nature of Sudoku has been widely studied, such as 
the brute force enumeration of solutions to $3 \times 3$ Sudoku \cite{Jarvis2006}, 
the sequential numbering of solutions \cite{Togami2006}，
density evolution analysis of belief-propagation based algorithms \cite{Jossy2014, Kitazono2015}，
Sudoku based nonlinear codes \cite{Jossy2015}, and 
the smallest number of clues required for a sSudoku puzzle to admit a unique solution \cite{McGuire2012, Reich2012}. 
\par
\begin{figure}[t]
  \begin{center}
    \includegraphics[width=0.60\linewidth,keepaspectratio]{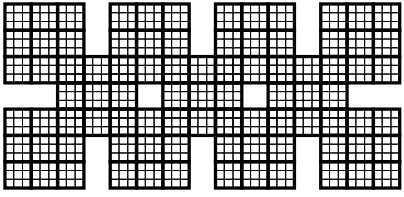} \\
    {\footnotesize (a) Shogun Sudoku grid.} \\[1em]
    \includegraphics[width=0.44\linewidth,keepaspectratio]{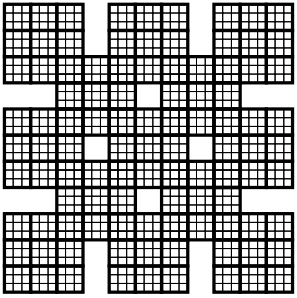} \\
    {\footnotesize (b) Sumo Sudoku grid.} \\[1em]
    \includegraphics[width=0.28\linewidth,keepaspectratio]{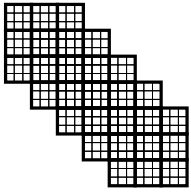} \\
    {\footnotesize (c) $\ell$-stage Stair Sudoku grid. The case of $\ell=5$.} \\[1em]
    \includegraphics[width=0.44\linewidth,keepaspectratio]{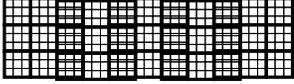} \\
    {\footnotesize (d) $\ell$-stage Belt Sudoku grid. The case of $\ell=5$.} 
    \caption{Examples of spatially coupled Sudokus.}
    \label{fig:scsudoku}
  \end{center}
\end{figure}

Here, we focus on spatially coupled Sudokus 
which consist of multiple ordinary Sudokus coupled by sharing some blocks. 
Examples of some spatially coupled sudokus are given in Figure \ref{fig:scsudoku}. 
This structure is closely related to spatially coupled low-density parity-check (LDPC) codes \cite{Kudekar2011,Yedla2014}．
The constraints for rows, columns, and blocks are individually applied to each Sudoku. 
In the study of Sudoku, the main concern is the number of solutions. 
When a Sudoku is regarded as an error correcting codes, it is used to evaluate the coding rate. 
\par
The main contribution of this paper is to evaluate upper bounds for various kinds of spatially coupled Sudokus. 
This paper is organised as follows. 
In Section \ref{sect:preliminaries}, we introduce some definitions and results from previous studies. 
In Section \ref{sect:mainresults}, the main results is presented. 
In Section \ref{sect:examples}, upper bounds for the number of solutions are evaluated 
for some typical spatially coupled Sudokus. 
A summary is provided in the final section.

\section{Preliminaries} \label{sect:preliminaries}
\par

\subsection{$n \times n$ Sudoku}
\par
Here, we consider an $n \times n$ Sudoku refers to composed of $n^2$ blocks of size $n \times n$. 
A {\it row band} refers to $n$ horizontally successive blocks 
and a {\it column band} refers to $n$ vertically successive blocks. 
In an $n \times n$ Sudoku, there are $n$ row bands and $n$ column bands.

\subsection{Permanent}
\par
To count solutions, we use the permanent of a matrix, according to Herzberg's analysis \cite{Herzberg2007}. 
For an $n \times n$ square matrix $A$ with $(i,j)$-th entry $a_{i,j}$, 
the permanent of $A$, which is denoted by $\mathrm{per} A$, is defined as 
\begin{align}
  \mathrm{per} A := \sum_{\sigma \in \mathcal{S}_n} \prod_{i=1}^n a_{i,\sigma(i)}, 
\end{align}
where $\mathcal{S}_n$ denotes the symmetric group for the $n$ symbols $\{1, \cdots, n\}$. 
Note that the permanent has a similar form to the determinant 
$\mathrm{det} A := \sum_{\sigma \in \mathcal{S}_n} \mathrm{sgn}(\sigma) \prod_{i=1}^n a_{i,\sigma(i)}$.  
\par
Let $A$ be an $n \times n$ $(0,1)$-matix with $r_i$ ones in row $i$, $i \in \{1, \cdots, n\}$. 
Then, $\mathrm{per} A$ is upperbounded as follows: 
\begin{align}
  \mathrm{per} A \le \prod_{i=1}^n r_i!^{1/r_i}. 
  \label{eq:per}
\end{align}
The proof of this inequality is given as Theorem 11.5 in \cite{VanLint1992}. 

\subsection{Previous Studies}
\par
Let $S(n)$ be the number of solutions to the $n \times n$ Sudoku. 
The following results have previously been obtained. 
\begin{theorem}[A part of Theorem 6 in \cite{Herzberg2007}]
  \label{Th:S(n)}
  The number of solutions to an $n \times n$ Sudoku is upperbounded by 
  \begin{align}
    S(n) \le 
    \prod_{i=1}^n 
    \biggl( \prod_{j=  1}^i \mu(i,j) \biggr) 
    \biggl( \prod_{j=i+1}^n \nu(j)   \biggr) 
    =: S_U(n), 
  \end{align}
  where 
  \begin{align}
    \mu(i,j) &:= [n^2-(i-1)n-(j-1)]!^{\frac{n^2}{n^2-(i-1)n-(j-1)}}, \\
    \nu(j) &:= [n^2-(j-1)n]!^{\frac{n^2}{n^2-(j-1)n}}. 
  \end{align}
  $S_U(n)$ is an upper bound for $S(n)$. 
  \qed
\end{theorem}
\begin{theorem}[Theorem 6 in \cite{Herzberg2007}]
  \label{Th:logS(n)}
  The upper bound $S_U(n)$ is given by 
  \begin{align}
    S_U(n) = n^{2n^4} \rme^{-2.5n^4 + O(n^3 \ln n)}, 
  \end{align}
  for sufficiently large $n$. 
  \qed
\end{theorem}

\subsection{Outline of the Proof of Theorems \ref{Th:S(n)} and \ref{Th:logS(n)}}
\par
Here, we briefly summarise Herzberg's analysis \cite{Herzberg2007}. 
The number of ways of completing the first row in the first row band is $n^2!$.
The number of ways of completing the second row in the first row band is evaluated 
by calculating the permanent of the following matrix. 
Let $A=(a_{i,j})$ be an $n^2 \times n^2$ $(0,1)$-matrix. 
The rows of $A$ parameterise the cells of the second row. 
The columns of $A$ parameterise the numbers from 1 to $n^2$. 
We set the $(i,j)$-th entry of $A$ to one if $j$ is a permissible value for cell $i$, and to zero otherwise. 
Then, $\mathrm{per}A$ gives the number of ways of choosing the set of distinct representatives. 
For instance, in the case that $n=3$ we can set 
\begin{align}
  A = \left(
  \begin{array}{ccc|ccc|ccc}
    0 & 0 & 0 & 1 & 1 & 1 & 1 & 1 & 1 \\
    0 & 0 & 0 & 1 & 1 & 1 & 1 & 1 & 1 \\
    0 & 0 & 0 & 1 & 1 & 1 & 1 & 1 & 1 \\ \hline
    1 & 1 & 1 & 0 & 0 & 0 & 1 & 1 & 1 \\
    1 & 1 & 1 & 0 & 0 & 0 & 1 & 1 & 1 \\
    1 & 1 & 1 & 0 & 0 & 0 & 1 & 1 & 1 \\ \hline
    1 & 1 & 1 & 1 & 1 & 1 & 0 & 0 & 0 \\
    1 & 1 & 1 & 1 & 1 & 1 & 0 & 0 & 0 \\
    1 & 1 & 1 & 1 & 1 & 1 & 0 & 0 & 0 \\
  \end{array}
  \right)
\end{align}
without loss of generality. 
Let $w_i(A) := \sum_{j=1}^{n^2} a_{i,j}$ be the weight of the $i$-th row of $A$. 
Setting $w_i(A)=n^2-n$ \,$\forall i$, 
the number of ways of completing the second row in the first row band 
can be evaluated as $(n^2-n)!^{n^2/(n^2-n)}$, by applying the inequality (\ref{eq:per}). 
Then we can obtain the number of ways of filling in the first row band consisting of $n$ rows as 
$\prod_{j=0}^{n-1} (n^2-jn)!^{n^2/(n^2-jn)}$. 
\par
Next, suppose that $(i-1)$ of the $n$ row bands have been completed. 
The number of possible entries for the first cell of the $i$-th row band is $n^2-(i-1)n$. 
Applying the column constraint, 
the number of possible entries for the $j$-th cell of the $i$-th row band becomes $n^2-(i-1)n-(j-1)$. 
On the other hand, by applying the block constraint this becomes $n^2-(j-1)n$. 
If $j \le i$, then $n^2-(i-1)n-(j-1) \le n^2-(j-1)n$. 
Using this property, we obtain the number of ways of filling in the $i$-th row band as 
$\{\prod_{j=1}^i [n^2-(i-1)n-(j-1)]!^{n^2/(n^2-(i-1)n-(j-1))} \}$
$\times \{\prod_{j=i+1}^n [n^2-(j-1)n]!^{n^2/(n^2-(j-1)n)} \}$. 
Thus, we arrive at Theorem \ref{Th:S(n)}. 
\par
Theorem \ref{Th:logS(n)} can be straightforwardly obtained 
by applying the Stirling's formula $\ln n! = n \ln n - n + O(\ln n)$ 
and some trivial inequalities to $\ln S_U(n)$.

\begin{figure}[t]
  \begin{center}
    \includegraphics[width=0.60\linewidth,keepaspectratio]{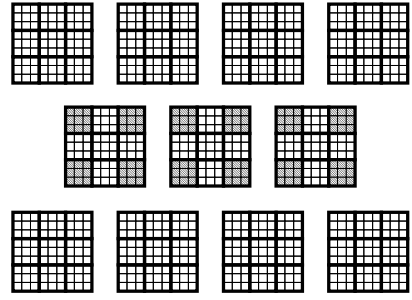} \\
    \caption{
      The Shogun Sudoku can be divided into eleven $3 \times 3$ Sudokus.
      A hatched block denotes a block in which all cells are already filled. 
    }
    \label{fig:division}
  \end{center}
\end{figure}

\section{Main Results} \label{sect:mainresults}
\par
In order to count the number of solutions to spatially coupled Sudokus, 
we first divide them into some isolated $n \times n$ Sudokus. 
Figure \ref{fig:division} illustrates an example of division. 
When spatially coupled Sudokus are divided into individual $n \times n$ Sudokus, 
some of the resulting isolated Sudoku can contain some blocks in which all cells are already filled. 
The set of positions of the filled blocks is not unique. 
We can choose a set of positions that minimises the upper bound.  
\par
\begin{figure}[t]
  \begin{center}
    \includegraphics[width=0.12\linewidth,keepaspectratio]{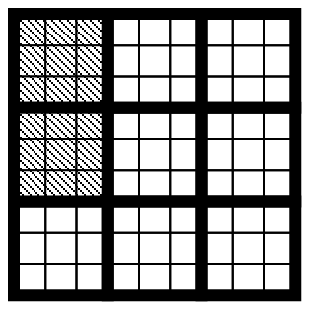} \\
    \caption{
      An $(n;c_1,c_2)$ partly filled Sudoku. The case of $(n;c_1,c_2)=(3;2,1)$. 
      All cells in the hatched $c_1 \times c_2$ rectangular shape blocks are filled. 
    }
    \label{fig:c1c2}
  \end{center}
  \vspace*{1em}
  \begin{center}
    \includegraphics[width=40mm]{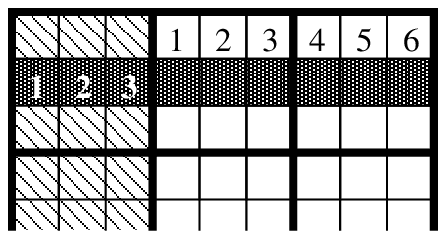}
    \caption{
      The case of $(n;c_1,c_2)=(3;2,1)$. 
      Only the first row band is shown. 
      To evaluate the upper bound, we consider the following situation without loss of generality. 
      The first cell in the second row of the first row band is 1. 
      The second cell is 2, 
      and the third cell is 3. 
      To ignore the effect of cell values that have been already filled, 
      we take only the block constraint into account. 
      The fourth to sixth cells are in $\{4,5,6,7,8,9\}$. 
      The seventh to nineth cells are in $\{1,2,3,7,8,9\}$. 
      For the other settings for $(n;c_1,c_2)$, one can consider a similar situation. 
    }
    \label{fig:sample}
  \end{center}
\end{figure}

Let $S(n;c_1,c_2)$ be the number of solutions to a partly filled $n \times n$ Sudoku 
whose $c_1 \times c_2$ blocks are filled as in Fig. \ref{fig:c1c2}. 
Such a Sudoku is referred to as an {\it $(n;c_1,c_2)$ partly filled Sudoku}. 
It should be noted that if any two row bands or any two column bands are interchanged, 
then the number of solutions to a Sudoku is unchanged, owing to its symmetric property. 
In detail, it can be explained as follows. 
For $U, V \subset I :=\{1, \cdots, n\}$, the Cartesian product $U \times V$ is called a {\it rectangle} in $I^2$. 
Let us consider a rectangle with the size of $|U|=c_1$ and $|V|=c_2$. 
The number of solutions to a partly filled Sudoku whose $(u,v)$-th blocks are already-filled for $\forall (u,v) \in U \times V$ 
depens only on the size of the rectangle, i.e., $S(n;,c_1,c_2)$ is independent of choice of $U$ and $V$.  
\par
\begin{theorem}[Upper bound for an $(n;c_1,c_2)$ partly filled Sudoku]
  \label{Th:m1}
  The number of solutions to an $(n;c_1,c_2)$ partly filled Sudoku is upperbounded by 
  \begin{align}
    S(n;c_1,c_2) \le 
    &
    \biggl\{ 
    \prod_{i=1}^{c_1}
    \biggl( \prod_{j=  1}^i \mu(i,j;c_2) \biggr) 
    \biggl( \prod_{j=i+1}^n \nu(j;c_2)   \biggr) 
    \biggr\} \notag\\
    & \times
    \biggl\{ 
    \prod_{i=c_1+1}^{n}
    \biggl( \prod_{j=  1}^i \mu(i,j;0) \biggr) 
    \biggl( \prod_{j=i+1}^n \nu(j;0)   \biggr) 
    \biggr\} \notag\\
    =:& 
    S_U(n;c_1,c_2), 
  \end{align}
  where 
  \begin{align}
    \mu(i,j;c_2) &:= [n^2-(i-1)n-(j-1)]!^{\frac{n^2-c_2n}{n^2-(i-1)n-(j-1)}}, \\
    \nu(j;c_2) &:= [n^2-(j-1)n]!^{\frac{n^2-c_2n}{n^2-(j-1)n}}. 
  \end{align}
  $S_U(n;c_1,c_2)$ is an upper bound for $S(n;c_1,c_2)$. 
  Note that $\mu(i,j;0)=\mu(i,j)$ and $\nu(j,0)=\nu(j)$. 
\end{theorem}
\par
{\it Proof:} 
We follow Herzberg's analysis \cite{Herzberg2007}. 
First, we consider the first $c_1$ row bands. 
The number of ways of completing the first row in the first row band is 
\[ (n^2-c_2n)!. \]
Considering an $n^2 \times n^2$ $(0,1)$-matrix $A$ which has the weight 
\begin{align*}
  w_i(A) = \left\{
  \begin{array}{ll}
    1, & \quad 1 \le i \le c_2n \\
    n^2, & \quad c_2n < i \le n^2 , 
  \end{array}
  \right. 
\end{align*}
and evaluating the upper bound (\ref{eq:per}) of $\mathrm{per} A$, 
this value is upperbounded by 
\[
  (n^2-c_2n)! 
  \le \mathrm{per} A
  \le n^2!^{\frac{n^2-c_2n}{n^2}}. 
\]
Although this makes the upper bound loose, it becomes easier to evaluate the upper bound for large $n$. 
The value $w_i(A)$ represents the number of possible digits that can be filled in the $i$-th cell of the first row in the first row band.
It should be noted that digits in already-filled blocks are not used to evaluate the upper bound 
since it is difficult to treat all possible cases that same digits in a block constraint appear in the corresponding row constraint; 
see Fig. \ref{fig:sample}. 
\par
The number of ways of completing the second row in the first row band is evaluated 
by calculating the permanent of the following $(0,1)$-matrix $A$. 
To evaluate the upper bound, we set 
\begin{align}
  w_i(A) = \left\{
  \begin{array}{ll}
    1, & \quad 1 \le i \le c_2n \\
    n^2 - n, & \quad c_2n < i \le n^2 , 
  \end{array}
  \right. 
\end{align}
without loss of generality. 
The number of ways of completing the second row in the first row band can be evaluated as 
\[ (n^2 - n)!^{\frac{n^2-c_2n}{n^2-n}}, \]
by applying the inequality (\ref{eq:per}). 
Considering a $(0,1)$-matrix $A$ which has the weight 
\begin{align*}
  w_i(A) = \left\{
  \begin{array}{ll}
    1, & \quad 1 \le i \le c_2n \\
    n^2 - 2n, & \quad c_2n < i \le n^2 ,
  \end{array}
  \right. 
\end{align*}
the number of ways of completing the third row in the first row band is evaluated as 
\[
  (n^2-2n)!^{\frac{n^2-c_2n}{n^2-2n}}. 
\]
Then, we obtain the number of ways of filling in the first row band consisting of $n$ rows as 
\begin{align}
  \prod_{j=1}^{n} [n^2-(j-1)n]!^{\frac{n^2-c_2n}{n^2-(j-1)n}}. 
\end{align}
\par
Next, suppose that $(i-1)$ of the $n$ row bands have been completed. 
The number of possible entries for the first cell of the $i$-th row band is $n^2-(i-1)n$. 
Applying the column constraint, 
the number of possible entries for the $j$-th cell of the $i$-th row band becomes $n^2-(i-1)n-(j-1)$. 
On the other hand, by applying the block constraint this becomes $n^2-(j-1)n$. 
If $j \le i$, then $n^2-(i-1)n-(j-1) \le n^2-(j-1)n$. 
Using this property, we obtain 
the number of ways of filling in the $i$-th row band as 
\begin{align}
  & \biggl( \prod_{j=1}^i [n^2-(i-1)n-(j-1)]!^{\frac{n^2-c_2n}{n^2-(i-1)n-(j-1)}} \biggr) \notag\\
  & \times \biggl( \prod_{j=i+1}^n [n^2-(j-1)n]!^{\frac{n^2-c_2n}{n^2-(j-1)n}} \biggr). 
\end{align}
\par
The remaining $n-c_1$ row bands can be treated in the same manner as in Theorem \ref{Th:S(n)}. 
We then arrive at Theorem \ref{Th:m1}. 
\qed
\par
By definition, it holds that $S(n;c_1,c_2) = S(n;c_2,c_1)$. 
Note that this upper bound takes into account all row constraints 
however, in terms of the column and the block constraint, only one of them is considered. 
Therefore, in regards to the result of Theorem \ref{Th:m1}, $S_U(n;c_1,c_2) \ne S_U(n;c_2,c_1)$ for $c_1 \ne c_2$ in general. 
\par
For large $n$, the exponent of the upper bound can be evaluated as follows. 
\begin{theorem}[Upper bound for a large $(n;c_1,c_2)$ partly filled Sudoku]
  \label{Th:m2}
  For $c_1=O(n^0)$ and $c_2=O(n^0)$, 
  the upper bound $S_U(n;c_1,c_2)$ is 
  \begin{align}
    S_U(n;c_1,c_2) = n^{2n^4} \rme^{-2.5n^4 + O(n^3 \ln n)}
  \end{align}
  for sufficiently large $n$. 
  Introducing parameters $d_1, d_2 \in [0,1]$, 
  let $c_1 = d_1 n = O(n) \in \mathbb{Z}$ 
  and $c_2 = d_2 n = O(n) \in \mathbb{Z}$. 
  The upper bound $S_U(n;d_1 n,d_2 n)$ is given by 
  \begin{align}
    S_U(n;d_1 n,d_2 n) = n^{\alpha(d_1,d_2)n^4} \rme^{\beta(d_1,d_2)n^4 + O(n^3 \ln n)}
  \end{align}
  for sufficienly large $n$, where
  \begin{align}
    \alpha(d_1,d_2) &:= 2 (1 - d_1 d_2), \\
    \beta(d_1,d_2) &:= - \frac 52
    + (1 - d_1) d_2 \ln (1 - d_1) 
    + d_1 d_2 
    + \frac{d_1^2 d_2}2. 
  \end{align}
  Here, $0 \ln 0 = 0$ by convention. 
\end{theorem}
\par
{\it Proof:} 
Theorem \ref{Th:m2} can be also obtained straightforwardly 
by applying the Stirling's formula $\ln n! = n \ln n - n + O(\ln n)$ 
and some trivial inequalities, i.e., 
$\sum_{j=0}^{i-1} \ln (n^2-(i-1)n-j) \le \sum_{j=0}^{i-1} \ln (n^2-(i-1)n)$ 
and $\ln(n-i) \le \ln (n-i+1)$, to $\ln S_U(n)$, which gives  
\begin{align}
  &
  \ln S_U(n;c_1,c_2) \notag \\
  &=
  n^2 \sum_{i=1}^{n}
  \biggl( \sum_{j=  1}^i \{ \ln (n^2-(i-1)n-j) - 1\} \notag\\
  &\quad + \sum_{j=i+1}^n \{ \ln (n^2-(j-1)n) -1 \} \biggr) \notag\\
  &\quad - c_2 n \sum_{i=1}^{c_1}
  \biggl( \sum_{j=  1}^i \{ \ln (n^2-(i-1)n-j) - 1\} \notag\\
  &\quad + \sum_{j=i+1}^n \{ \ln (n^2-(j-1)n) -1 \} \biggr) 
  + O(n^2 \ln n) \notag\\
  &\le 
  2n^4 - \frac 52 n^2 -c_2 n^3 \ln \frac{n}{n-c_1} - c_1 c_2 n^2 \ln n(n-c_1) \notag\\
  &\quad - c_1(2n-c_1)\frac{c_2}2 n + 2 c_1 c_2 n^2
  + O(n^3 \ln n). 
\end{align}
\qed
For $c_1=O(n^0)$ and $c_2=O(n^0)$, 
the effects of filled blocks are neglected. 
Namely, the upeer bound depends on neither $c_1$ nor $c_2$. 
When $d_1$ or $d_2$ is equal to zero, 
it holds that $\alpha(0,d_2) = \alpha(d_1,0) = 2$ and $\beta (0,d_2) = \beta (d_1,0) = -5/2$.

\section{Examples} \label{sect:examples}
\par
We apply our result to some typical spatially coupled Sudokus and obtain the following results. 
The coding rate can be defined as an error correcting code as follows. 
\begin{definition}[Coding Rate]
  The coding rate $R$ of a Sudoku can be described as an error correcting code by 
  \begin{align}
    R = \frac {\log_{n^2} S}{C}, 
  \end{align}
  where $S$ and $C$ denote the number of solutions to the Sudoku 
  and the number of cells in the Sudoku grid, respectively. 
  \qed
\end{definition}
Let $R(n;c_1,c_2)$ be the coding rate of an $(n;c_1,c_2)$ partly filled sudoku. 
\begin{example}[$2 \times 2$ Sudoku]
  The number of solutons to a $2 \times 2$ Sudoku $S(2;0,0)$ can be easily obtained 
  by brute force or simple counting. 
  We compare it to the upperbound $S_U(2;0,0)$ as follows: 
  \begin{align}
    S(2;0,0) = 288 \le S_U(2;0,0) = 384. 
  \end{align}
  The coding rate of the is 
  \begin{align}
    R(2;0,0) = \frac{\log_4 S(2;0,0)}{16} \approx 0.2553. 
  \end{align}
  The number of solutons to some partly filled Sudokus can be easily obtained as 
  \begin{align}
    & S(2;1,1) = 12 \le S_U(2;1,1) = 39, \\
    & S(2;1,2) = 4 = S_U(2;1,2). 
  \end{align}
  by brute force or simple counting. 
  Note that in the case of $(n;c_1,c_2)=(2;1,2)$, 
  the number of solutions $S(2;1,2)$ depends on a pattern of cell values that have already filled 
  and can take 2 or 4. 
  \qed
\end{example}
\begin{example}[$3 \times 3$ Sudoku]
  Felgenhauer and Jarvis obtained the number of solutons 
  to a $3 \times 3$ Sudoku $S(3;0,0)$ by brute force \cite{Jarvis2006}. 
  We compare the result to the upperbound $S_U(3;0,0)$: 
  \begin{align}
    S(3;0,0) 
    &= 6, 670, 903, 752, 021, 072, 936, 960 \notag \\
    &\approx 6.6709 \times 10^{21} \notag\\
    &\le S_U(3;0,0) \approx 1.7071 \times 10^{26}. 
  \end{align}
  The coding rate becomes 
  \begin{align}
    R(3;0,0)  
    = \frac{\log_9 S(3;0,0)}{81}
    \approx 0.2823. 
  \end{align}
  \qed
\end{example}
\begin{example}[The Shogun Sudoku]
  The number of solutions to the Shogun Sudoku $S^{shogun}$ is upperbounded by 
  \begin{align}
    S^{shogun}
    &\le S(3;0,0)^8 S(3;2,2)^3 \notag\\
    &\le S(3;0,0)^8 S_U(3;2,2)^3 =: S_U^{shogun} \notag\\
    &\approx (6.6709 \times 10^{21})^8 \times (1.5976 \times 10^{11})^3 \notag\\
    &\approx 1.5993 \times 10^{208}. 
  \end{align}
  The upper bound of the coding rate $R_U^{shogun}$ can be evaluated as 
  \begin{align}
    R^{shogun}
    &\le \frac{\log_9 S_U^{shogun}}{783} =: R_U^{shogun} \notag\\
    &= \frac{\log_9 1.5993 \times 10^{208}}{783} \approx 0.2786.
  \end{align}
  \qed
\end{example}
\begin{example}[The Sumo Sudoku]
  The number of solutions to the Sumo Sudoku $S^{sumo}$ is upperbounded by 
  \begin{align}
    S^{sumo}
    &\le S(3;0,0)^9 S(3;2,2)^4 \notag\\
    &\le S(3;0,0)^9 S_U(3;2,2)^4 =: S_U^{sumo} \notag\\
    &\approx (6.6709 \times 10^{21})^9 \times (1.5976 \times 10^{11})^4 \notag\\
    &\approx 1.7045 \times 10^{241}. 
  \end{align}
  The upper bound of the coding rate $R_U^{sumo}$ is 
  \begin{align}
    R^{sumo}
    &\le \frac{\log_9 S_{U,sumo}}{909} =: R_U^{sumo} \notag\\
    &= \frac{\log_9 1.7045 \times 10^{241}}{909} \approx 0.2781.
  \end{align}
  \qed
\end{example}
\begin{example}[The $\ell$-stage Stair Sudoku]
  The number of solutions to the $\ell$-stage stair Sudoku $S^{stair}(\ell)$ is upperbounded by 
  \begin{align}
    S^{stair}(\ell)
    &\le 
    S(3;0,0) S(3;2,2)^{\ell-1} \notag\\
    &\le 
    S(3;0,0) S_U(3;2,2)^{\ell-1} =: S_U^{stair}(\ell) \notag\\
    &\approx 
    10^{21.8241 + 11.2034 (\ell-1)}. 
  \end{align}
  The upper bound of the coding rate $R_U^{stair}(\ell)$ is
  \begin{align}
    R^{stair}(\ell)
    &= \frac{\log_9 S_{U,stair}(\ell)}{81+54(\ell-1)} =: R_U^{stair}(\ell) \notag\\
    &\approx \frac{22.8706 + 11.7407 (\ell-1)}{81+45(\ell-1)}. 
  \end{align}
  In the large $\ell$ limit, it becomes $\lim_{\ell \to \infty }R_U^{stair}(\ell) \approx 0.2609$. \\
  \qed
\end{example}
\begin{example}[The $\ell$-stage Belt Sudoku]
  The upper bound for 
  the number of solutions to the $\ell$-stage belt Sudoku $S^{belt}(\ell)$ is given by 
  \begin{align}
    S^{belt}(\ell)
    &\le 
    S(3;0,0) S(3;1,3)^{\ell-1} \notag\\
    &\le 
    S(3;0,0) S_U(3;1,3)^{\ell-1} =: S_U^{belt}(\ell) \notag\\
    &\approx 
    10^{21.8241 + 14.0520 (\ell-1)}. 
  \end{align}
  Note that $S_U(3;1,3) \le S_U(3;3,1)$. 
  The upper bound of the coding rate $R_U^{belt}(\ell)$ is
  \begin{align}
    R^{belt}(\ell)
    &\le \frac{\log_9 S_U^{belt}(\ell)}{81+54(\ell-1)} =: R_U^{belt}(\ell) \notag\\
    &\approx \frac{22.8706 + 14.7258 (\ell-1)}{81+54(\ell-1)}. 
  \end{align}
  We then have $\lim_{\ell \to \infty }R_U^{belt}(\ell) \approx 0.2727$. 
  \qed
\end{example}

\section{Summary} \label{sect:summary}
\par
We have evaluated the upper bounds for the number of solutions to spatially coupled Sudokus 
such as the Shogun Sudoku, the Sumo Sudoku, $\ell$-stage stair Sudoku, and $\ell$-stage belt Sudoku. 
Brute force enumeration and 
evaluation of the lower bound and tighter upper bound will be the focus of our future studies.

\section*{Acknowledgments}
\par
This work was partially supported by 
a Grant-in-Aid for Scientific Research 
(B) Nos. 16K12496 \& 25289114, 
(C) No. 25330264, and
for Challenging Exploratory Research No. 16K12496
from the Ministry of Education, Culture, Sports, Science and Technology (MEXT) of Japan.



\begin{thebibliography}{9}


\bibitem{Jarvis2006}
B. Felgenhauer, F. Jarvis, 
``Mathematics of Sudoku I,'' 
{\it Mathematical Spectrum}, vol. 39, pp. 15--22, 2006. 


\bibitem{Togami2006}
S. Togami, 
bachelor's degree thesis, 
Tokyo institute of Technology, 2006. 

\bibitem{Jossy2014}
C. Atkins and J. Sayir, 
``Density Evolution for SUDOKU Codes on the Erasure Channel,'' 
{\it Proc. of ISTC2014}, pp. 233--237, Aug. 2014. 

\bibitem{Kitazono2015}
T. Kitazono and K. Mimura, 
``Coding Theoretical Evaluation for Spatially Coupled Sudokus,'' 
{\it Proc. of HISS2015}, Nov. 2015 (in Japanese). 

\bibitem{Jossy2015}
J. Sayir and J. Sarwar, 
``An investigation of SUDOKU-inspired non-linear codes with local constraints,''  
{\it Proc. of ISIT 2015}, pp. 1921--1925, Jun. 2015.

\bibitem{McGuire2012}
G. McGuire, B. Tugemann, and G. Civario, 
``There is no 16-Clue Sudoku: Solving the Sudoku Minimum Number of Clues Problem via Hitting Set Enumeration,'' 
{\it Experimental Mathematics}, vol. 23, no. 2, pp. 190--217, 2014; arXiv http://arxiv.org/abs/1201.0749, Jan. 2012. 

\bibitem{Reich2012}
E. S. Reich, 
``Mathematician claims breakthrough in Sudoku puzzle,'' 
{\it Nature News}, doi:10.1038/nature.2012.9751, Jan. 2012

\bibitem{Kudekar2011}
S. Kudekar, T. J. Richardson, and R. L. Urbanke, 
``Threshold Saturation via Spatial Coupling: Why Convolutional LDPC Ensembles Perform So Well over the BEC,'' 
{\it IEEE Trans. Info. Theory}, vol. 57, no. 2, pp. 803--834, Feb. 2011. 

\bibitem{Yedla2014}
A. Yedla, Y.-Y. Jian, P. S. Nguyen, and H. D. Pfister, 
``A Simple Proof of Maxwell Saturation for Coupled Scalar Recursions,'' 
{\it IEEE Trans. Info. Theory}, vol. 60, no. 11, pp. 6943--6965, Aug. 2014. 


\bibitem{Herzberg2007}
A. M. Herzberg and M. Ram Murty, 
``Sudoku squares and chromatic polynomials,'' 
{\it Notices of the AMS}, 54, pp. 708–717, 2007.

\bibitem{VanLint1992}
J. H. Van Lint and R. M. Wilson, 
``A Course in Combinatorics,'' 
{\it Cambridge University Press}, 1992. 



\end{thebibliography}
\end{document}